\documentclass[final,11pt]{article}
\usepackage{amsmath}
\usepackage{amsfonts}
\usepackage{amsthm}

\def\R{\mathbb{R}}

\begin{document}

\title{A symmetry result for strictly convex domains
}

\author{A. G. Ramm\\
\small Mathematics Department\\
\small Kansas State University, Manhattan, KS 66506, USA\\
\small email: ramm@math.ksu.edu}
\date{}
\maketitle

\begin{abstract}

Assume that $D\subset \R^2$ is a strictly convex domain
with $C^2-$smooth boundary.

{\bf Theorem.} {\em If $\int_De^{ix}y^ndxdy=0$ for
all sufficiently large $n$, then $D$ is a disc.}

\end{abstract}



\section{Introduction} \label{Introduction}

We assume throughout that $D\subset \R^2$ is a {\em strictly convex} domain
 and its boundary $S$ is  $C^2-$smooth.
Suppose that
\begin{equation}\label{eq:2}
 \int_D e^{ix}y^ndxdy=0, \qquad n=0,1,2,......
\end{equation}
Our result is stated as Theorem 1.

{\bf Theorem 1.} {\em If $D$ is strictly convex bounded domain in $\mathbb{R}^2$
and (\ref{eq:2}) holds, then $D$ is a disc.}

This result the author obtained while studying the Pompeiu problem, see, for example,
Chapter 11 in the book \cite{R470}. The result of Theorem 1 can also be established
if the following is assumed in place of equation (\ref{eq:2}):
\begin{equation}\label{eq:3}
 \int_D e^{iy}x^ndxdy=0, \qquad n=0,1,2,......
\end{equation}
 This follows from the proof of Theorem 1.

\section{Proof of Theorem 1.}

 Let $\ell$ be an arbitrary unit vector, $L_1$ be the support line to $D$ (at the point $s_1\in S$)
parallel to $\ell$,
and $L_2$ be  the support line to $D$  (at the point $q_1\in S$)  parallel to $L_1$, where $q_1=q_1(s_1)$.
 Since $D$ is strictly convex, one can introduce the equations $y=f(x)$ and  $y=g(x)$
of the boundary $S$ between the support points $s_1$ and $q_1$.
For definiteness and without loss of generality
 let us assume that the orthogonal projection of the point $s_1$ onto the line  $L_1$ lies not
lower than the projection of the point $q_1$ onto $L_1$, and let the $x-$axis pass through $s_1$
and be orthogonal to $L_1$. The graph of $f$ is located above
the graph of $g$.  Since $S$ is strictly convex the function $f$ has a unique point of maximum
$x_1$, where $x_1\in (a,b)$, and $f(x_1)>f(x)$ for $x\in [a,b]$, $f(x_1)>0$ and $f^{''}(x_1)<0$.
Here $a$ and $b$ are the $x-$coordinates of the points  $q_1$ and $s_1$, $a<b$.
Let us denote by $s$ the value of the natural parameter (arc length on $S$) corresponding to the
maximum point of $f$, that is, to the point $x_1$. The function $g$ has a unique point of  minimum
 $x_2$,  $x_2\in (a,b)$, $g(x)>g(x_2)$,  $g(x_2)<0$ and  $g^{''}(x_2)>0$.
  From the strict convexity of $S$ it follows that these maximum and minimum are
non-degenerate, that is, $f^{''}(x_1)\neq 0,$ and $g^{''}(x_2)\neq 0.$
Denote by $q$ the value of the natural parameter corresponding to the minimum point of $g$.
Let us write   formula (\ref{eq:2}) as
 \begin{equation}\label{eq:6}
\int_D e^{ix}y^ndxdy=\int_a^b e^{ix} \frac {f^{n+1}(x)-g^{n+1}(x)}{n+1} dx= 0, \qquad n=0,1,2,.......
\end{equation}
The factor $n+1$ in the denominator can be canceled because the integral in (\ref{eq:6}) equals to zero.
We want to take $n\to \infty$ and use the Laplace method for evaluating the main term of
the asymptotic of the integral. Let us recall this known result, the formula
 for the asymptotic of the integral
 \begin{equation}\label{eq:7'}
 F(\lambda):=\int_a^b \phi(x) e^{\lambda S(x)}dx= \Big(\frac{2\pi}{\lambda |S^{''}(\xi)|}\Big)^{1/2} \phi(\xi)e^{\lambda S(\xi)}
 \Big (1+o(1)\Big),\qquad \lambda\to \infty,
 \end{equation}
see, for example, \cite{B}. In this formula
 $\xi\in (a,b)$  is a unique point of a non-degenerate maximum of a real-valued twice continuously
 differentiable function $S(x)$ on $[a,b]$,  $S^{''}(\xi)< 0$,
 and $\phi$ is a continuous function on $[a,b]$, possibly complex-valued.
We apply this formula with
$$S(x)=\ln |f|,\quad \lambda:=2m:=n+1\to \infty,\quad \phi=e^{ix},$$
and take $n=2m-1$ to ensure that $n+1=2m$ is an even number, so that
$f^{2m}$ and $g^{2m}$ are positive,  and
$\ln f^{2m}$ and $\ln g^{2m}$ are well defined.  The point $x_2$ of minimum
of $g$ becomes a point of local maximum of the function   $g^{2m}$.
Note that $|(\ln |f|)^{''}|=
\frac {|f^{''}(x_1)|}{|f(x_1)|}$ at the point $x_1$ where $f'(x_1)=0$, $f(x_1)>0$ and $f^{''}(x_1)<0$.

Taking the above into consideration, one obtains from (\ref{eq:6}) the following asymptotic formula:
%\begin{equation}\label{eq:7}
\begin{eqnarray}
\int_D e^{ix}y^ndxdy=\Big[e^{ix_1+2m\ln |f(x_1)|}\Big(\frac{\pi|f(x_1)|}{m|f^{''}(x_1)|}\Big)^{1/2}- \notag\\
e^{ix_2+2m\ln |g(x_2)|}\Big(\frac{\pi|g(x_2)|}{m|g^{''}(x_2)|}\Big)^{1/2}\Big]\Big(1+o(1)\Big)= 0, \quad n\to \infty,
\end{eqnarray}
where $2m=n+1$, $x_1\in (a,b)$ and $x_2\in (a,b)$.
%\end{equation}
 It follows from the above formula  that the expression in the brackets, that is, the main term of the asymptotic,
  must vanish for all sufficiently large $m$. This implies that $f(x_1)=|f(x_1)|=|g(x_2)|$ and
  $|f^{''}(x_1)|=g^{''}(x_2)=|g^{''}(x_2)|$, because $f(x_1)>0$, $g(x_2)<0$, $f^{''}(x_1)<0$
  and $g^{''}(x_2)>0$. It also follows from formula  (5) that $e^{ix_1}=e^{ix_2}$. This implies $x_1=x_2+2\pi p$,
  where $p$ is an integer. The integer $p$ does not depend on $s$ because $p$ is locally
  continuous and cannot have jumps.
 Thus,
 \begin{equation}\label{eq:8}
x_1-x_2:=2\pi p; \quad |f(x_1)|= |g(x_2)|; \quad |f^{''}(x_1)|=|g^{''}(x_2)|.
\end{equation}
We prove in Lemma 2 (see below) that $p=0$. Another proof of this is given in the Remark 1 below the
proof of Lemma 2.

Consider the support lines $L_3$ at the point $s$ and $L_4$ at the point $q$, where $L_3$ and $L_4$
are orthogonal to $\ell$.  Denote by $L=L(s)$
the distance between $L_3$ and $L_4$, that is, the width
of $D$ in the direction parallel to $\ell$. Note that $L=f(x_1)-g(x_2)>0$, and
 \begin{equation}\label{eq:9}
L=(r(s)-r(q), \ell),
\end{equation}
 where  $r=r(s)$ is the radius vector (position vector)
corresponding to the point on $S$ which is defined by the parameter $s$. This point will be called point $s$.
 The same letter $s$
is used for the point $s\in S$ and for the corresponding natural parameter.
Let $R=R(s)$ denote the radius of curvature of the curve $S$ at the point $s$ and
let $\kappa=\kappa (s)$ denote the curvature of $S$ at this point.  Then one has
 \begin{equation}\label{eq:9'}
R^{-1}=\kappa=|f^{''}(x_1)|,
\end{equation}
because $\kappa=|f^{''}(x_1)|[1+|f'(x_1)|^2]^{-\frac{3}{2}} $
and $f'(x_1)=0$ since $x_1$ is a point of maximum of $f$.

 From   (\ref{eq:8}) we will derive that
\begin{equation}\label{eq:10}
L(s)=2R(s), \qquad \forall s\in S.
\end{equation}
It will be proved in  Lemma 2, see below, that  equation (\ref{eq:10})
implies that  $D$ is a disc.
Thus, the conclusion of Theorem 1 will be established.

  We denoted by $r=r(s)$  the equation of $S$, where $s$ is the natural parameter  on $S$
and  $r$ is the radius vector of the point on $S$, corresponding to $s$.
One has $r'(s)=t$, where $t=t(s)$ is a unit vector tangential to $S$
 at the point $s$. We have chosen $s$ so that $t(s)$
is orthogonal to $\ell$. Since $\ell$ is arbitrary, the point $s\in S$ is arbitrary.
The point $q\in S$, $q=q(s)$,  is uniquely determined by the requirement that
$t(q)=-t(s)$, because $S$ is strictly convex.
% Vector $r(q)-r(s)$ is
%orthogonal to the vectors $t(s)$ and $t(q)$, and
One has $(r(s)-r(q), \ell)=L$, where $L=L(s)$
is the width of $D$ in the direction parallel to $\ell$.
Since $ r'(s)=t(s)$, the first formula  (\ref{eq:8}) implies
\begin{equation}\label{eq:14}
(r(q)-r(s), r'(s))=2\pi p,\quad (r(q)-r(s), r'(q))=- 2\pi p, \quad \forall s\in S.
\end{equation}
Differentiate the first equation  (\ref{eq:14}) with respect to $s$
and get
\begin{equation}\label{eq:15}
(r'(q)\frac{dq}{ds}-r'(s), r'(s))+(r(q)-r(s), r^{''}(s)) =0, \quad \forall s\in S.
\end{equation}
 Note that $r'(s)=t(s)=-t(q)=-r'(q)$ and  $ r^{''}(s)=\kappa(s) \nu (s)$,
 where $\nu(s)$ is the unit normal to $S$ (at the point corresponding to  $s$)
directed into $D$, and $(r(s)-r(q), \ell)=L(s)=(r(q)-r(s), \nu(s))$,
because $\nu(s)$ is directed along $-\ell$.  Consequently, it follows from  (\ref{eq:15}) that
\begin{equation}\label{eq:16}
-\frac{dq}{ds}-1+\kappa(s)L(s) =0, \quad \forall s\in S.
\end{equation}
One has $L(s)=L(q)$, and it follows from formulas (\ref{eq:8}) that  $\kappa(s)=\kappa(q)$.

Differentiate the second equation  (\ref{eq:14}) with respect to $q$ and get
\begin{equation}\label{eq:16'}
(t(q)-t(s)\frac{ds}{dq}, t(q))+(r(q)-r(s), r^{''}(q)) =0.
\end{equation}
Note that $t(q)=-t(s)$ and $ r^{''}(q)=\kappa(q)\nu(q)$, where $\nu(q)=-\nu(s)$
 because $L_3$ is parallel to $L_4$. Consequently, equation (\ref{eq:16'}) implies
 \begin{equation}\label{eq:17}
\frac{ds}{dq}+1-\kappa(s)L(s) =0, \quad \forall s\in S.
\end{equation}
Compare   (\ref{eq:16}) and  (\ref{eq:17}) and get $\frac{ds}{dq}=\frac{dq}{ds}$.
Thus, $\Big(\frac{dq}{ds}\Big)^2=1$. Since $\frac{dq}{ds}>0$, it follows that
\begin{equation}\label{eq:18}
\frac{ds}{dq}=\frac{dq}{ds}=1, \qquad \forall s\in S.
\end{equation}
Therefore, equation  (\ref{eq:16})  implies
\begin{equation}\label{eq:19}
\kappa (s)L(s)=2 \qquad \forall s\in S.
\end{equation}
{\em Let us derive from (\ref{eq:19}) that $D$ is a disc.}

Recall that
$s$ is the natural parameter on $S$, $L(s)$ is the width of $D$
at the point $s$ (that is the distance between two parallel supporting lines to $S$
one of which passes through the point $s$) and $\kappa(s)$ is the curvature of $S$ at
the point $s$.

{\bf Lemma 2.} {\it Assume that $D$ is strictly convex domain with a smooth boundary $S$.
 If  equation  (\ref{eq:19}) holds, then $D$ is a disc.}

{\it Proof of Lemma 2.}
Denote by $K$ the maximal
disc inscribed in the strictly convex domain $D$, and by $r$ the radius of $K$.
If there are no points of $S$ outside $K$, then $D$ is a disc and we are done.
If $S$
contains points outside $K$, let $x\in S$ be such a point.
Consider the line  $\tilde{L}$ passing through the center of $K$ and
through the point $x\in S$, $x\not\in K$. Let $L^{'}$ be the support line to $S$
 orthogonal to the line $\tilde{L}$ and tangent to $S$ at a point $x'$, $x'\not\in K$.
 Denote the radius of curvature of $S$ at the point $x'$
 by $\rho$. One has $\rho\le r$, because $K$ is the maximal disc inscribed in $D$.
The width $L$ of $D$ at the point $x'$ in the direction
of the line  $\tilde{L}$  is greater than $2r$ because $x'\not\in K$. One has $L>2r$
and $L=2\rho \le 2r$.  This is a contradiction. It
 proves that $D=K$.  Thus,  $D$ is a disc, and, consequently, the parameter $p$ in formula
(\ref{eq:8}) is equal to zero.
Lemma 2 is proved. \hfill $\Box$

Thus, Theorem 1 is proved. \hfill $\Box$

{\bf Remark 1.} Let us give another proof that $p=0$, where $p$ is defined in formula (\ref{eq:8}).
One has $L(s)=(r(q)-r(s),\nu(s))$. Differentiate this equation with respect to $s$ and
get
\begin{equation}\label{eq:20}
L'(s)=-(r(q)-r(s), \kappa (s) t(s))+ (r'(q)\frac{dq}{ds}-r'(s),\nu(s)),
\end{equation}
where $t(s)$ is
the unit vector tangential to $S$ at the point $s$. Here the known formula $\nu(s)'=-\kappa(s)t(s)$
was used.  The second term in the formula (\ref{eq:20})
vanishes since $r'(s)$ and $r'(q)$ are orthogonal to $\nu$. Thus,   $L'(s)=-2\pi p \kappa(s)$.
Since $D$ is strictly convex, one has inequality $min_{s\in S}\kappa (s)\ge \kappa_0>0$, where
$\kappa_0>0$ is a constant. The function $L(s)$ must be periodic, with the period
 equal to the arc length of $S$. The differential equation  $L'(s)=-2\pi p \kappa(s)$ does not have periodic
solutions unless $p=0$. Therefore, $p=0$. \hfill $\Box$

\end{document}